\documentclass[11pt,letterpaper]{amsart}

\usepackage{amsmath, amssymb,latexsym, amsthm, amsfonts, amscd}
\usepackage{verbatim}
\newtheorem{thm}{Theorem}

\pdfoutput=1

\usepackage{graphics}

\usepackage{citesort}
\usepackage{amssymb,amsmath}
\usepackage[mathscr]{eucal}

\newtheorem{lemma}{Lemma}
\normalsize

\title[Knots and Links Arising From Recombination]{A Topological Characterization Of Knots and Links Arising From Site-Specific Recombination}

\author{Dorothy Buck and Erica Flapan}
\date \today

\address{Department of Mathematics and Centre for Bioinformatics, Imperial College of London, UK}
\email{d.buck@imperial.ac.uk}
\address{Department of Mathematics, Pomona College,
Claremont, CA 91711, USA}
\email{eflapan@pomona.edu}

\newcommand{\nibf}{\noindent \textbf}

\newcommand{\beq}{\begin{equation}}
\newcommand{\eeq}{\end{equation}}
\newcommand{\bea}{\begin{eqnarray}}
\newcommand{\eea}{\end{eqnarray}}

\begin{document}

\bigskip
\begin{abstract}
We develop a topological model of knots and links arising from a
single (or multiple processive) round(s) of recombination starting
with an unknot, unlink, or $(2,m)$-torus knot or link substrate.
We show that all knotted or linked products fall into a single
family, and prove that the size of this family grows linearly with
the cube of the minimum number of crossings. Additionally, we
prove that the only possible products of an unknot substrate are
either clasp knots and links or $(2,m)$-torus knots and links.
Finally, in the (common) case of $(2,m)$-torus knot or link
substrates whose products have minimal crossing number $m+1$, we
prove that the types of products are tightly prescribed, and use
this to examine previously uncharacterized experimental data.
\end{abstract}

\maketitle

\section{Introduction}

Molecular biologists are interested in DNA knots and links,
because they have been implicated in a number of cellular
processes. The axis of DNA molecules can become knotted or linked
as a result of many reactions, including replication and
recombination. The wide variety of DNA knots and links observed
has made separating and characterizing these molecules a critical
issue. Experimentally, this is most conclusively accomplished via
electron microscopy \cite{KrasStas}.  However, this is a laborious
and difficult process. Thus topological techniques, such as those
presented here, can aid experimentalists in characterizing DNA
knots and links by restricting the types of knots or links that
can arise in a particular context.

This work focuses on one such DNA knotting process,
\textit{site-specific recombination}, mediated by a
protein, known as a \textit{site-specific recombinase}.
Site-specific recombination is important because of
its key role in a wide variety of biological processes. (See
\cite{BFbio} or \textit{e.g.} \cite{MobDNA} for more information).
In addition, pharmaceutical and agricultural industries have become increasingly involved in genetically modifying organisms or testing whether a mutation in a particular gene leads to a disease.  As a result, these industries are now interested in site-specific recombinases as tools for precisely manipulating DNA (\textit{e.g.} \cite{Feil}).

Site-specific recombination roughly has three stages. Two recombinase molecules first bind to each
of two specific sites on one or
two molecules of covalently closed circular DNA (known as the
\textit{substrate}) and then bring them close
together. We shall refer to these DNA sites as the \textit{crossover
sites}. Next, the sites are cleaved, exchanged and resealed. The
precise nature of this intermediary step is determined by which of
the two recombinase subfamilies the particular protein belongs to
(see Assumption 3 below for more details). And finally, the
rearranged DNA, called the \textit{product}, is released.
Understanding precisely which knots and links arise during
site-specific recombination can help understand the details of
this process (\textit{e.g.} \cite{FlpAnti}).

In this paper we begin by developing a model that predicts all possible
knots and links which can arise as products of a single round of
recombination, or multiple rounds of (processive) recombination,
starting with substrate(s) consisting of an unknot, an unlink, or
a $(2,m)$-torus knot or link (denoted by $T(2,m)$). This model rigorously develops and extends ideas that we originally sketched in  \cite{OCAMI}. Of all knots and links, we have chosen to focus on
$T(2,m)$, because $T(2,m)$ are the most commonly occurring knots and links in DNA .  Our model rests on three assumptions that we
justify biologically in \cite{BFbio}. Building on
these assumptions, we use knot theoretic techniques to prove
that all products fall within a single family, illustrated in
Figure \ref{productfamily}.   We then prove that the number of product knots and links predicted by our model grows linearly with the cube of the minimal crossing number.
We further prove that the product
knot or link type is tightly prescribed when the substrate is $T(2,m)$ and the product has minimal crossing number $m+1$.  Finally, we apply this new result
to previously uncharacterized experimental data.

This paper complements earlier work by Sumners, Ernst, Cozzarelli
and Spengler \cite{SECS}, which used the tangle model \cite{ES1}
and several biologically reasonable assumptions to solve tangle
equations. They then determined which 4-plat knots and links arise
as a result of (possibly processive) site-specific recombination
on the unknot for the serine subfamily of recombinases (see just before Assumption 3 for a discussion of the two subfamilies). For the particular case of
the recombinase Gin, they considered the knots $3_1,4_1,5_2$ or $6_1$ as well as unknotted
substrates. Our paper goes further in several ways. In
addition to allowing an unknotted substrate for a generic recombinase, we
allow substrates that are unlinks with one site on each component,
as well as any $T(2,m)$. Also, our
assumptions are based exclusively on the biology of the
recombination process. In particular, we do not assume the tangle
model holds or that all products must be 4-plats. Allowing products which are not 4-plats is important because recombination has been seen to produce knots and
links which are connected sums (see \cite{BFbio}).

\begin{figure}[htpb]
\includegraphics{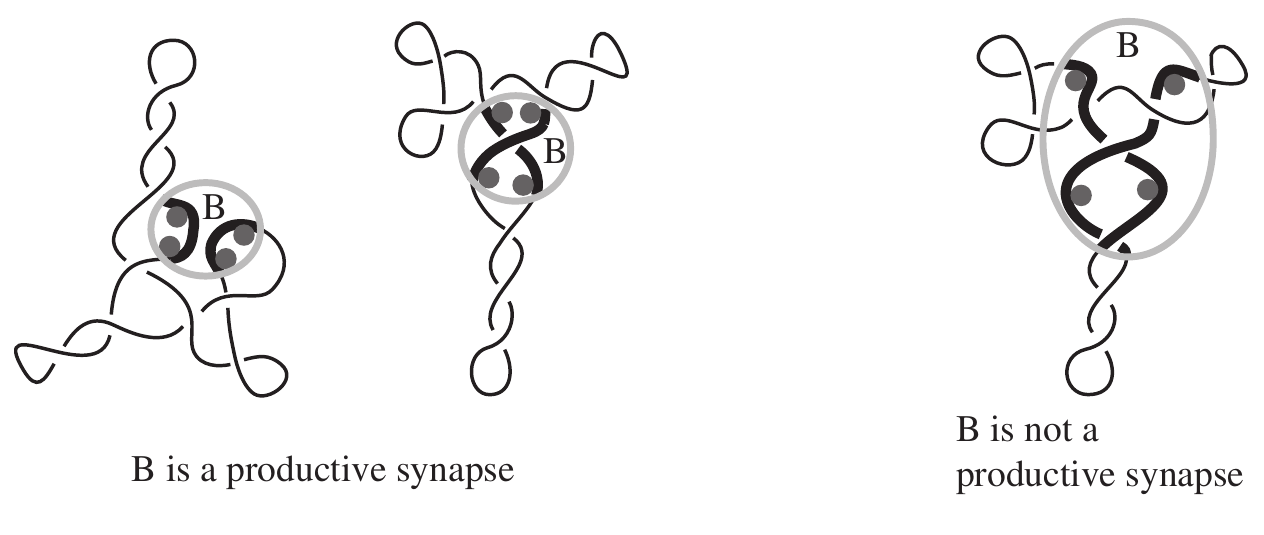}
\caption{The two examples on the left have a productive synapse
and the one on the right does not. The crossover sites are
highlighted.} \label{productivesynapse}
\end{figure}

We will use the following terminology and notation. Let $J$ denote
a substrate which is either an unknot, an unlink, or $T(2,m)$
(illustrated in Figure \ref{BJCD}). We use the term
\textit{recombinase complex, B,} to refer to the convex hull of
the four bound recombinase molecules together with the two
crossover sites, and use the term {\it recombinase-DNA complex} to
refer to $B$ together with the substrate $J$. If the recombinase
complex meets the substrate in precisely the two crossover sites
then we say the recombinase complex is a {\it productive synapse}.
In Figure \ref{productivesynapse}, we illustrate two examples
where the recombinase complex $B$ is a productive synapse, and one
where $B$ is not. Finally, we let $C=\mathrm{cl}(\mathbb{R}^3-B)$,
and consider $B\cap J$ and $C\cap J$ separately.

The structure of the paper is as follows. In Section 2, we state
our three assumptions about the recombinase-DNA complex, and use
our assumptions to determine the pre-recombinant and
post-recombinant forms of $B\cap J$. In Section 3, we characterize
the forms of $C\cap J$ for each of our substrates. In Section 4,
we glue each of the post-recombinant forms of $B\cap J$ to each
form of $C\cap J$ to determine all possible knotted or linked
products predicted by our model. Finally in Section 5, we bound the size of this product family, and further limit product type in some special cases, by combining our model with results on minimal crossing number.

\medskip

\section{Our assumptions and $B\cap J$}

\subsection{The three assumptions}
We make the following three assumptions about the recombinase-DNA complex, which we state in both biological and mathematical terms.  In \cite{BFbio}
we provide experimental evidence showing that each of these assumptions is biologically reasonable.

\medskip

\nibf{(Biological) Assumption 1:} The recombinase complex is
a productive synapse, and there is a projection of the crossover
sites which has at most one crossing between the sites and no
crossings within a single site.

\smallskip

This is equivalent to:

\smallskip

\noindent{\bf (Mathematical) Assumption 1:}  $B\cap J$ consists of two arcs and there is a projection of
$B\cap J$ which has at most one crossing between the two arcs, and
no crossings within a single arc.

\begin{figure}[h]
\includegraphics{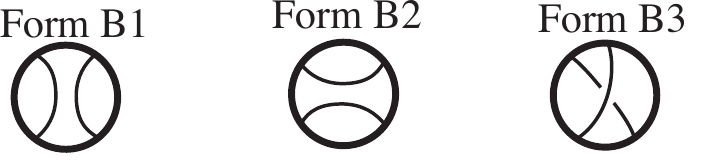}
\caption{ We fix a projection of $J$ so that $B\cap J$ has one of these forms.}
\label{BcapJ}
\end{figure}
As a result of this assumption, we now fix a projection of $J$ such that $B\cap J$ has one of the forms illustrated in Figure \ref{BcapJ}.

\medskip

\noindent{\bf(Biological) Assumption 2:} The productive synapse
does not pierce through a supercoil or a branch point in a
nontrivial way. Also, no persistent knots are trapped in the
branches of the DNA on the outside of the productive synapse.

\medskip

\begin{figure}[htpb]
\includegraphics{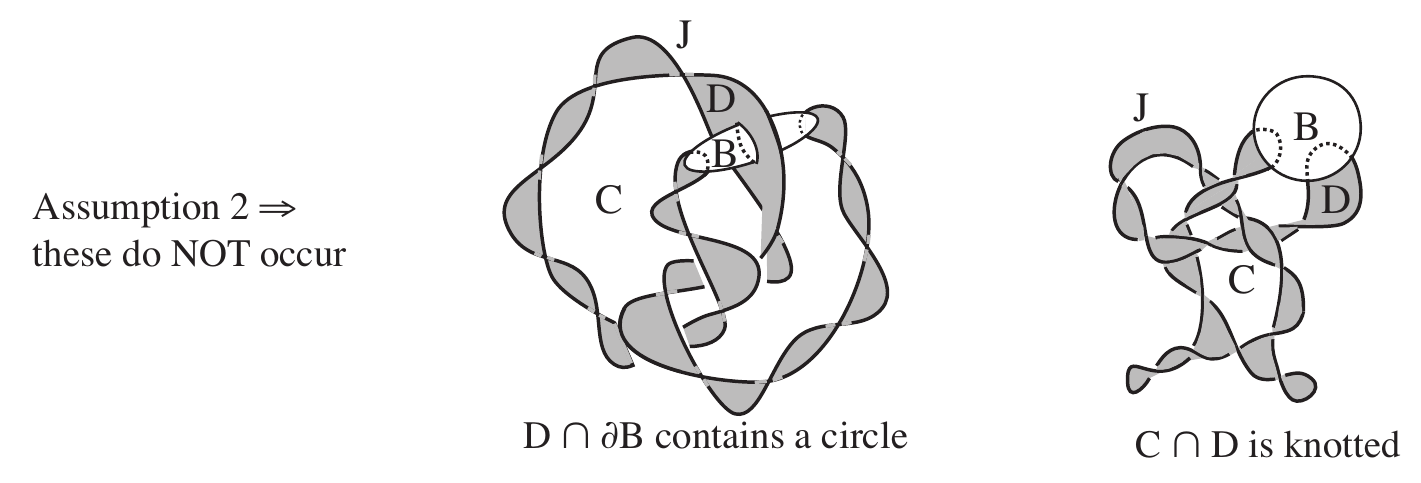}
\caption{On the left, the productive synapse pierces through a
supercoil in a nontrivial way. On the right, a knot is trapped in
the branches on the outside of $B$} \label{piercePaper}
\end{figure}

Assumption~2 implies that the
recombinase-DNA complex cannot resemble either of the illustrations
in Figure \ref{piercePaper}.

\begin{figure}[h]
\includegraphics{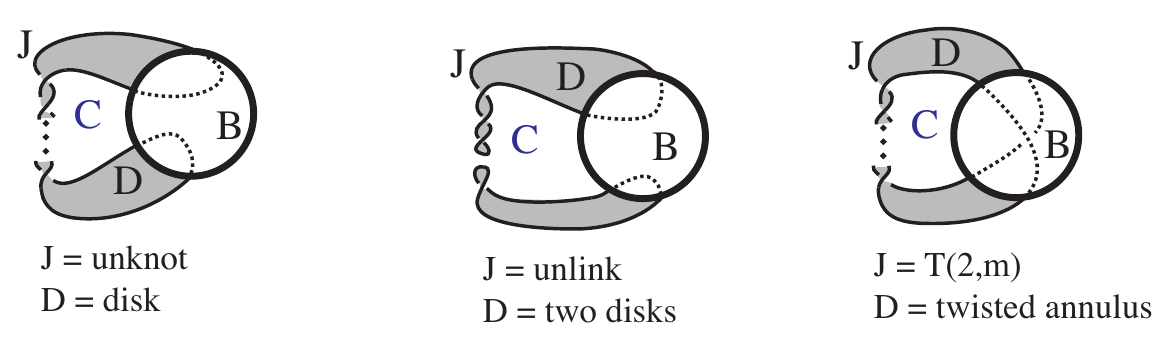}
\caption{ Examples of different substrates $J$ and a spanning surface $D$ bounded by $J$.}
\label{BJCD}
\end{figure}
In order to restate Assumption 2 mathematically, we first introduce some terminology.  We shall use the term {\it
spanning surface} to refer to a surface $D$, bounded by $J$, such that $D$ is topologically equivalent to a disk, two
disjoint disks, or a twisted annulus when $J$ is an unknot, unlink, or $T(2,m)$, respectively.  Figure \ref{BJCD} gives some examples of
the relationship between a spanning surface $D$ and the productive synapse $B$.  Observe that in each
of the illustrations of
Figure \ref{BJCD}, $D\cap \partial B$ consists of two arcs.  By Assumption~1, $B$ contains precisely two arcs of $J=\partial D$.
Hence $\partial B$ meets $J$ in precisely four points. It follows that the intersection of any spanning surface for $J$ with $\partial B$ contains
exactly two arcs.  What we mean by  $B$ {\it does not pierce through a
supercoil or a branch point in a nontrivial way} is that $B$ does not
pierce the interior of every spanning surface for $J$ (as in the
left illustration in Figure \ref{piercePaper}). In general, a spanning
surface $D$ is pierced by $B$ if and only if $D\cap \partial B$ contains at least one circle in addition to the required two arcs.
For example, in the diagram on
the left in Figure \ref{piercePaper}, no matter how the spanning
surface $D$ is chosen, $D\cap
\partial B$ contains at least one circle as well as two arcs.

\begin{figure}[h]
\includegraphics{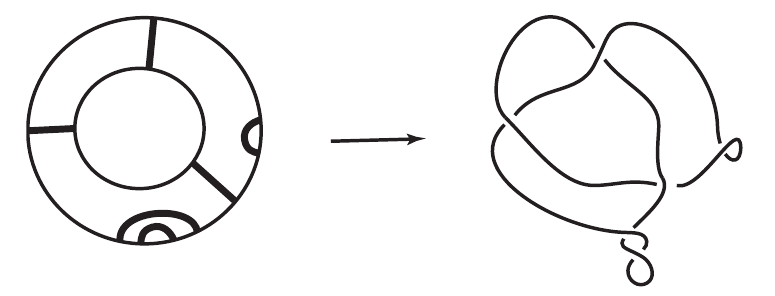}
\caption{We obtain a {\it planar surface with twists} by replacing a neighborhood of each arc by a half-twisted band.} \label{planar}
\end{figure}

Next, we explain what we mean by no {\it persistent knots are trapped in the branches outside of $B$}. Consider a planar surface together with a finite
number of arcs whose endpoints are on the boundary of the surface (see the illustration on the left in Figure \ref{planar}).
We can obtain a surface bounded by a knot or link by replacing a neighborhood of each arc in the original surface
by a half-twisted band and removing the top and bottom ends of the band. Figure \ref{planar} illustrates how such a surface
can be obtained from an annulus together with a collection of
arcs defining the twists. Any surface obtained from a planar surface in
this way is said to be a {\it planar surface with twists}.

Suppose that $D$ is a spanning
surface for $J$.  We say that $D\cap C$ is {\it unknotted rel $\partial B$}, if
there is an ambient isotopy of $C$ pointwise fixing $\partial B$ which takes $D\cap C$ to a planar
surface with twists, where the endpoints of the arcs defining the twists are disjoint from $\partial B$. For example, $D\cap C$ is unknotted rel $\partial B$ for each
of the spanning surfaces in Figure \ref{BJCD}. This is not the case for
the surfaces $D\cap C$ in Figure \ref{piercePaper}.

We now restate Assumption~2 mathematically as follows.

\medskip

\noindent{\bf (Mathematical) Assumption 2:} $J$ has a spanning
surface $D$ such that $D\cap \partial B$ consists of two arcs and
$D\cap C$ is unknotted rel $\partial B$.

\medskip

Site-specific recombinases fall into two families -- the serine
and tyrosine recombinases. Assumption~3 addresses the mechanism of
recombination according to which subfamily the recombinase is in.
While the overall reactions of the two families of recombinases
are the same, they differ in their precise mechanism of cutting
and rejoining DNA at the crossover sites. We explain more of the
biological details in \cite{BFbio}.
\medskip

\noindent {\bf (Biological) Assumption 3:} {\it Serine}
recombinases perform recombination via the ``subunit exchange
mechanism." This mechanism involves making two simultaneous
(double-stranded) breaks in the sites, rotating opposites sites
together by $180^{\circ}$ within the productive synapse and
resealing opposite partners. In processive recombination, each
recombination event is identical. After recombination mediated by
a {\it tyrosine} recombinase, there is a projection of the
crossover sites which has at most one crossing.

\medskip

The mathematical restatement of Assumption 3 is almost identical to the biological statement.
\medskip

\noindent {\bf (Mathematical) Assumption 3:} {\it Serine}
recombinases cut each of the crossover sites and add a crossing
within $B$ between the cut arcs on different sites, then reseal.
In processive recombination, all recombination events are
identical. After recombination mediated by a {\it tyrosine}
recombinase, there is a projection of the crossover sites which
has at most one crossing.
\medskip

\subsection{The forms of $B\cap J$}
\begin{figure}[h]
\includegraphics{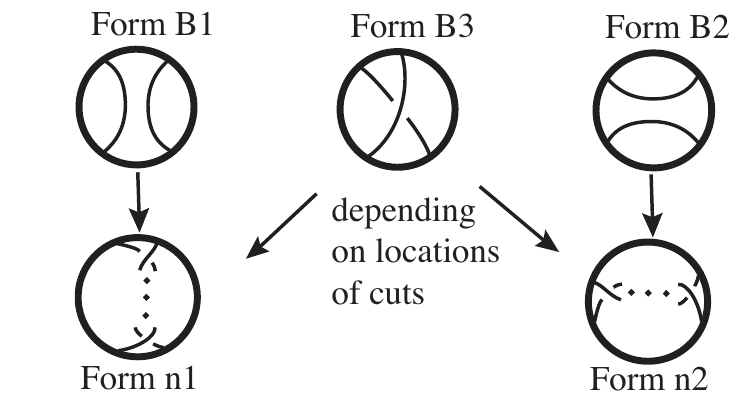}
\caption{After $n$ recombination events with serine recombinases, each pre-recombinant form of $B\cap J$ leads to the corresponding post-recombinant form.} \label{nBJPaper}
\end{figure}

As a result of Assumption~1, we fixed a projection of $J$ such
that $B\cap J$ has Form B1, B2, or B3 (illustrated in Figure
\ref{BcapJ}). It follows from Assumption 3 that after $n$
recombination events with serine recombinases, we have added a row
of $n$ identical crossings. Thus after $n$ recombination events our fixed projection of $B\cap J$ is isotopic
fixing $\partial B$ to one of the forms illustrated in Figure
\ref{nBJPaper} (where the actual crossings can be positive, negative, or
zero).

\begin{figure}[h]
\includegraphics{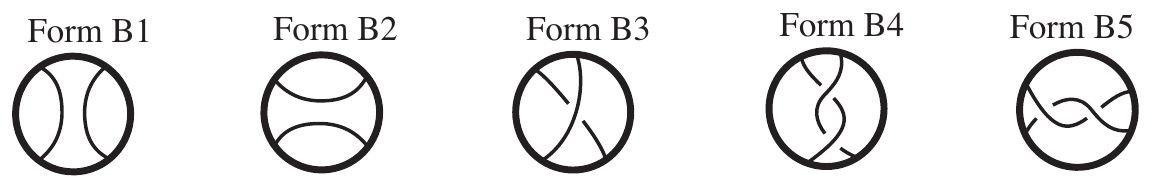}
\caption{After recombination with tyrosine recombinases, the post-recombinant projection of $B\cap J$ has one of these forms.} \label{tyroBJ1}
\end{figure}

Also for tyrosine recombinases, we know from Assumption~3 that after recombination there exists a projection of $B\cap J$ with at most one crossing.  We are working with the projection of $J$ which we fixed as a result of Assumption 1, and we cannot be sure that this particular projection of $B\cap J$ will have at most one crossing.  However our projection must be ambient isotopic, fixing $\partial B$, to one of the forms illustrated in Figure \ref{tyroBJ1}.  So without loss of generality we will assume that the post-recombinant projection of $B\cap J$ has one of these forms.

\medskip

\section { The possible forms of $C\cap J$}

Using Assumption 2, we now prove the following Lemma.

\begin{figure}[htpb]
\includegraphics{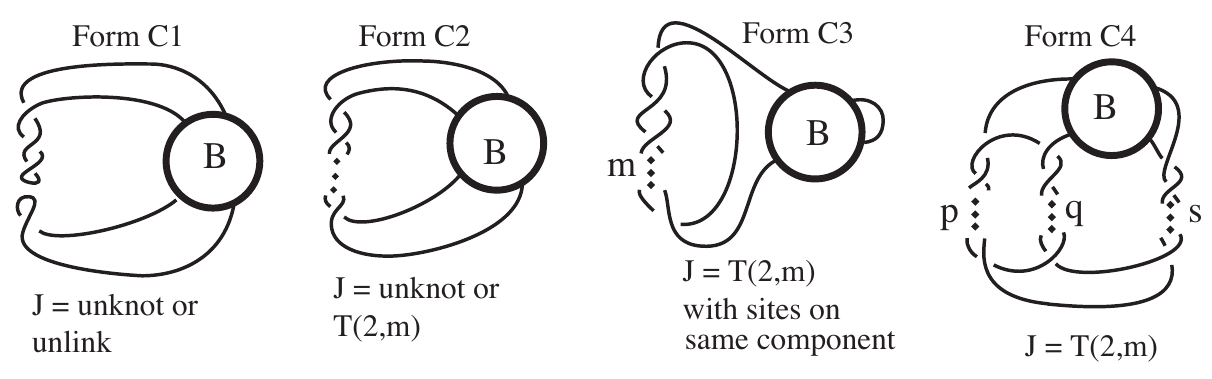}
\caption{$C\cap J$ has a projection with one of these forms.}
\label{newforms}
\end{figure}

 \begin{lemma} \label{L:CcapJ} Suppose that Assumptions 1, 2, and 3 hold for a particular
recombinase-DNA complex where the substrate $J$ is an unknot, unlink, or a $T(2,m)$ knot or link.  Then $C\cap J$ has a projection with one of the forms illustrated in Figure \ref{newforms} where $p+q=m$.  Furthermore, if $C\cap J$ has Form C4, then $B\cap J$ must have Form B1 in Figure \ref{tyroBJ1}.
\end{lemma}

 \begin{proof}

 We consider separate cases according
to the knot or link type of $J$.

\medskip

 \noindent {\bf Case 1:}  $J$ is the unknot.

In this case, by Assumption~2, we can choose a spanning surface
$D$ which is a disk such that $D\cap \partial B$ is two arcs and $D\cap C$ is unknotted rel $\partial B$. Since $D$ is a disk, the two arcs of $\partial B \cap D$ separate $D$ such that one of $B\cap D$ and $C\cap D$
is a strip and the other is a pair of disjoint disks.  Furthermore, if $C\cap D$ is a strip it is not knotted. Thus, $C\cap D$ is either a pair of disjoint disks
or an unknotted twisted strip. It follows that $C\cap J$ is ambient isotopic, pointwise fixing $\partial B$, to Form C1 or Form C2.

 \medskip

\noindent {\bf Case 2:}   $J$ is the unlink

In this case, we assume that one site is on each component of $J$
(or else the substrate was actually an unknot). Thus by
Assumption~2, we can choose a spanning surface $D$ which is a pair of disjoint disks
such that $ \partial B$ meets each disk of $D$ in a single arc.
Hence, $B\cap D$ and $C\cap D$ are each a pair of disjoint
disks. It follows that $C\cap J$ is ambient isotopic, pointwise fixing $\partial B$, to Form C2.

\medskip

\noindent {\bf Case 3:}   $J=T(2,m)$

In this case, by Assumption~2, we can choose a spanning surface
$D$ to be a twisted annulus such that $D\cap \partial B$ is two arcs and $D\cap C$ is unknotted rel $\partial B$. We see as follows that there are several ways the arcs of $D\cap \partial B$ can lie in $D$.

Any circle in $\mathbb{R}^3$ must cross a sphere an
even number of times (possibly zero). In particular, the circle
$A$ representing the core of the twisted annulus $D$ must cross
$\partial B$ an even number of times. Each point where $A$ crosses
$\partial B$ is contained in $D\cap\partial B$. Since the total
number of points in $A\cap \partial B$ is even and $D\cap\partial
B$ consists of two arcs, either $A$ must intersect each of these two
arcs an even number of times, or $A$ must intersect each of the two
arcs an odd number of times. If each arc of $D\cap\partial
B$ intersects the core $A$ an odd number of times, then each of these arcs cuts $D$ into a strip. Hence the two
arcs of $D\cap \partial B$ together cut $D$ into a pair of strips. If each arc of $D\cap\partial
B$ intersects the core $A$ an even number of times, then each arc cuts off a disk from $D$.
In this case,  either the two arcs cut off disjoint disks in $D$, or one of
the disks is contained inside of the other. In this latter case,
the two arcs form the edges of a strip in $D$, on one side of which is a disk and on the other side of which is a twisted annulus. The
three forms of $D\cap \partial B$ are illustrated on the top of
Figure \ref{torusforms}.  Note that the illustration on the right may have one, rather than two, rows of twists.  Since $B\cap J$ contains at most one crossing, the component of $D$ with almost all of the twists of $T(2,m)$ must be contained in $C$.

\begin{figure}[h]
\includegraphics{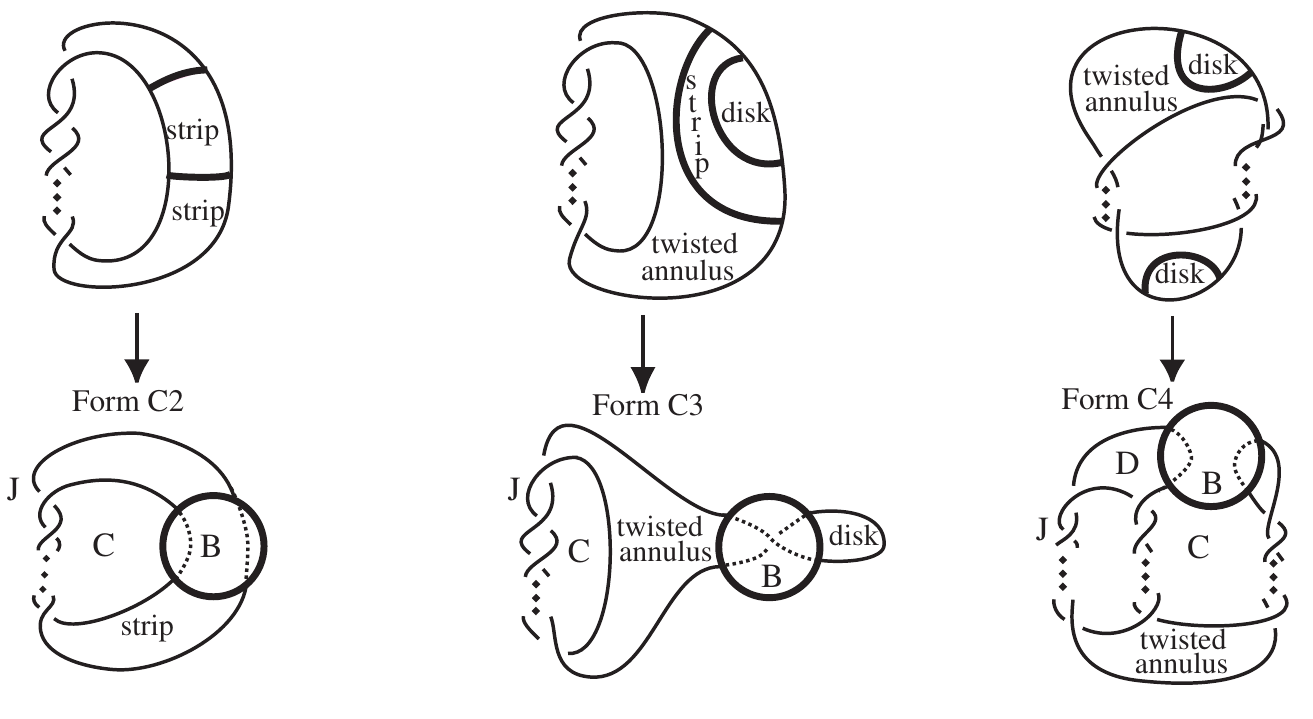}
\caption{These are the forms of $C\cap D$ when $J=T(2,m)$.}
 \label{torusforms}
\end{figure}

Since $C\cap D$ is unknotted rel $\partial B$, the abstract forms illustrated on the top of Figure \ref{torusforms} yield the corresponding forms of $C\cap D$ which are illustrated in the bottom of Figure \ref{torusforms} up to isotopy fixing $\partial B$.  Observe that when $C\cap J$ has Form C4, the projection of $B\cap J$ must have Form B1 as illustrated.  Also, in Form C3, while there may be twists to the left of $B$, they are topologically insignificant, since they can be removed by rotating $D\cap C$ by some multiple of $\pi$.  Similarly, in Form C4, any twists which had occurred to the left of $B$ can be removed and added to the row of twists at the right by rotating $D\cap C$ by some multiple of $\pi$.  These rotations can occur while pointwise fixing $B$.

Thus the four forms of $C\cap D$ illustrated in Figure \ref{newforms} are the only ones possible.
\end{proof}
\medskip

\section{Product knots and links predicted by our model}

\begin{figure}[h]
\includegraphics{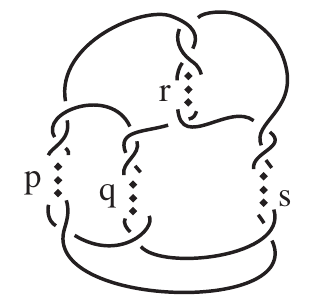}
\caption{We show that all knotted and linked products are in this family.}
\label{productfamily}
\end{figure}

In this section, we suppose that the substrate is an unknot, an
unlink, or $T(2,m)$ and that all three of our assumptions hold for
a particular recombinase-DNA complex. Then we prove Theorems
\ref{T:tyrosine} and \ref{T:serine}, which characterize all
possible knotted or linked products brought about by tyrosine recombinases and serine recombinases respectively. If the substrate is an unknot or unlink we will also
show that all nontrivial products are in the torus link family $T(2,n)$ or the clasp link family
$C(r,s)$ (i.e., consisting of one row of $r$ crossings and a
non-adjacent row of $s$ crossings). Note that $C(r,\pm2)$ is the well known family of
{\it twist} knots and links. If the substrate is $T(2,m)$, then
all products are in the family of knots and links illustrated in
Figure \ref{productfamily}.

\begin{figure}[htpb]
\includegraphics{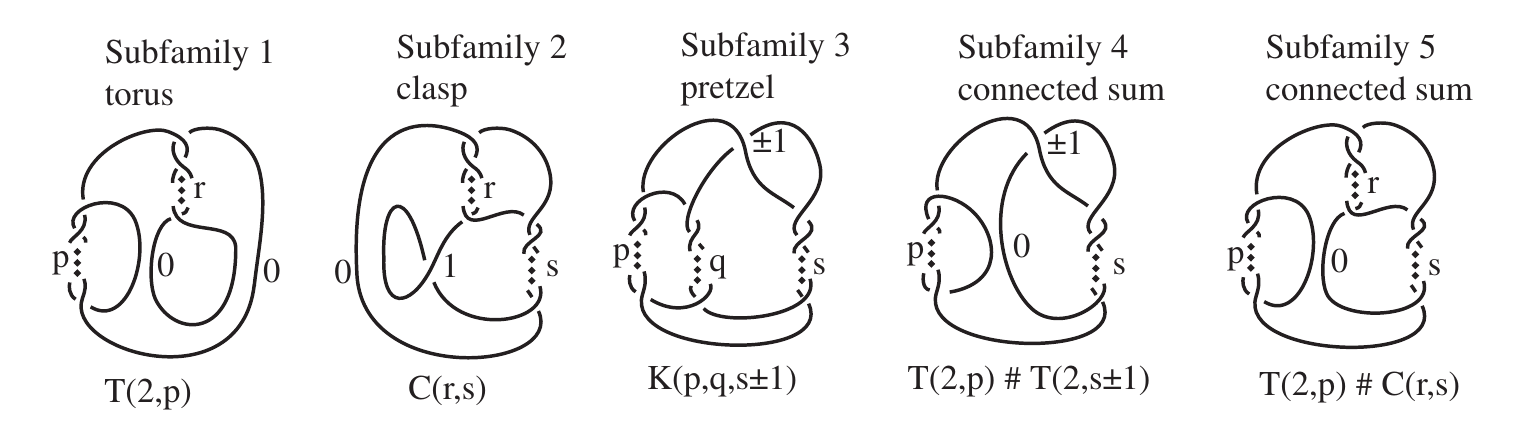}
\caption{These subfamilies are contained in the family illustrated in Figure \ref{productfamily}.}
\label{family5}
\end{figure}

Observe that  in Figure \ref{productfamily}, $p$, $q$,
$r$, and $s$ can be positive, negative, or zero.  Furthermore, by letting $p$, $q$, $r$, and/or $s$ equal 0 or 1 as appropriate, we obtain the five
subfamilies illustrated in Figure \ref{family5}. Subfamily 3 is the family of pretzel knots or links $K(p,q,r)$ with three non-andjacent rows
containing $p$ crossings, $q$ crossings, and $r$ crossings respectively.  Observe that Subfamily 4 is a connected sum.  However, if $q=0$, $r=1$, and $s=-1$, then it is a $T(2,p)$ together with an unlinked trivial component.

\begin{thm}\label{T:tyrosine}
Suppose that Assumptions 1, 2, and 3 hold for a particular
tyrosine recombinase-DNA complex. If the substrate is an unknot
then the only nontrivial products are $T(2,n)$ or $C(2,n)$. If
the substrate is an unlink, then the only nontrivial product is a
Hopf link, $T(2,2)$. If the substrate is $T(2,m)$,
then all of the products are contained in the family illustrated
in Figure \ref{productfamily}.
\end{thm}

\begin{figure}[h]
\includegraphics{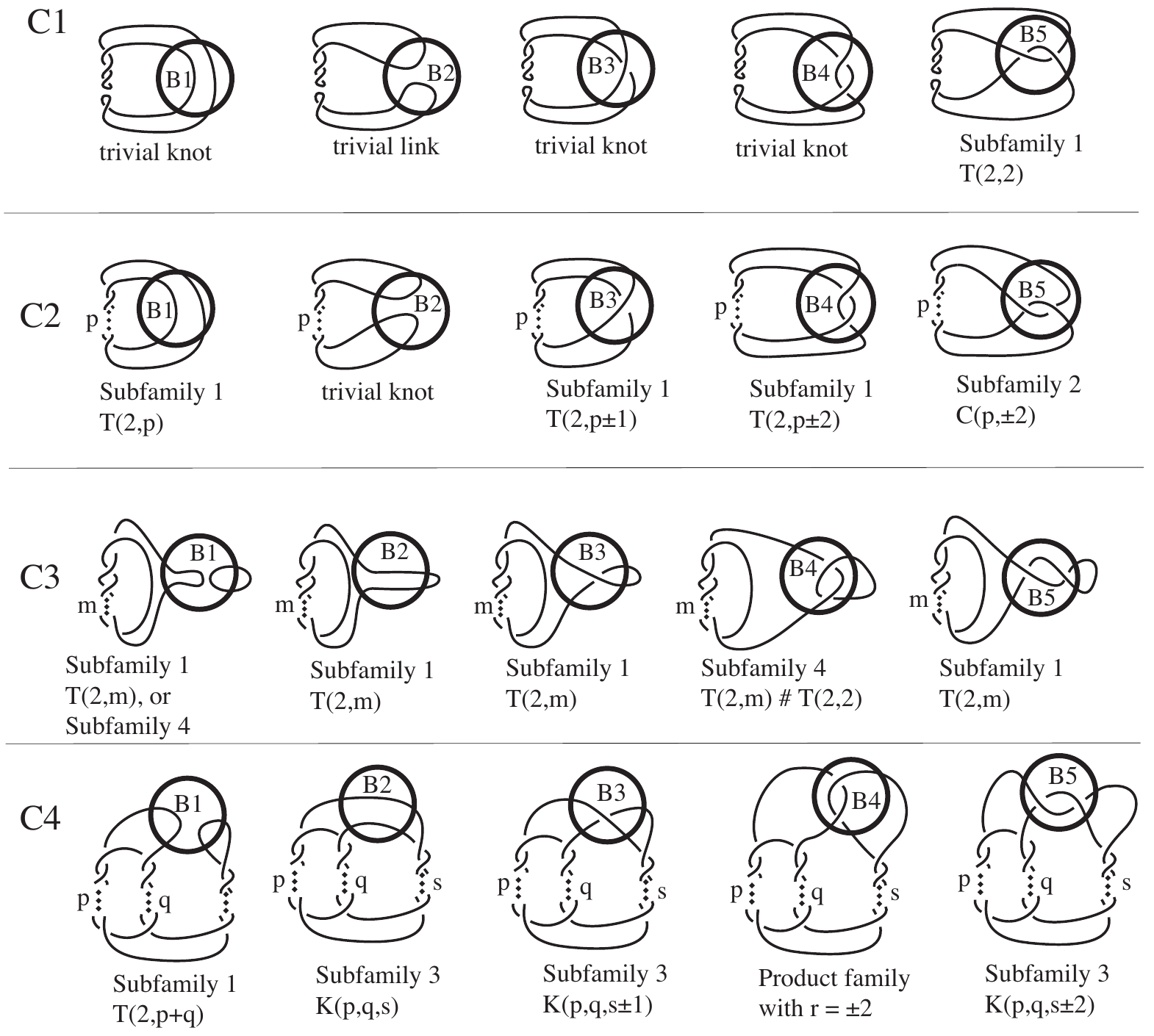}
\caption{Products of recombination with tyrosine recombinases.}
\label{tyrosine}
\end{figure}

\begin{proof}

We saw that as a result of Assumption~3, after recombination with tyrosine recombinases, the fixed projection of $B\cap J$ is ambient isotopic fixing $\partial B$ to one of the five forms illustrated in
Figure \ref{tyroBJ1}.  Also, by Lemma \ref{L:CcapJ}, $C\cap J$ is ambient isotopic fixing $\partial B$ to one of the four forms illustrated in Figure
\ref{newforms}.  For each of the four forms of $C\cap J$, the products of recombination with tyrosine recombinases are obtained by
replacing $B$ with each of the five post-recombinant forms of $B\cap J$ in Figure \ref{tyroBJ1}. The resulting products
are illustrated in Figure \ref{tyrosine}. Recall that if $J$ is the unlink then
$C\cap J$ must have Form C1, and if $J$ is an unknot then $C\cap J$ must have Form C1 or C2.  The theorem follows from Figure \ref{tyrosine}.

\end{proof}

\begin{figure}[h]
\includegraphics{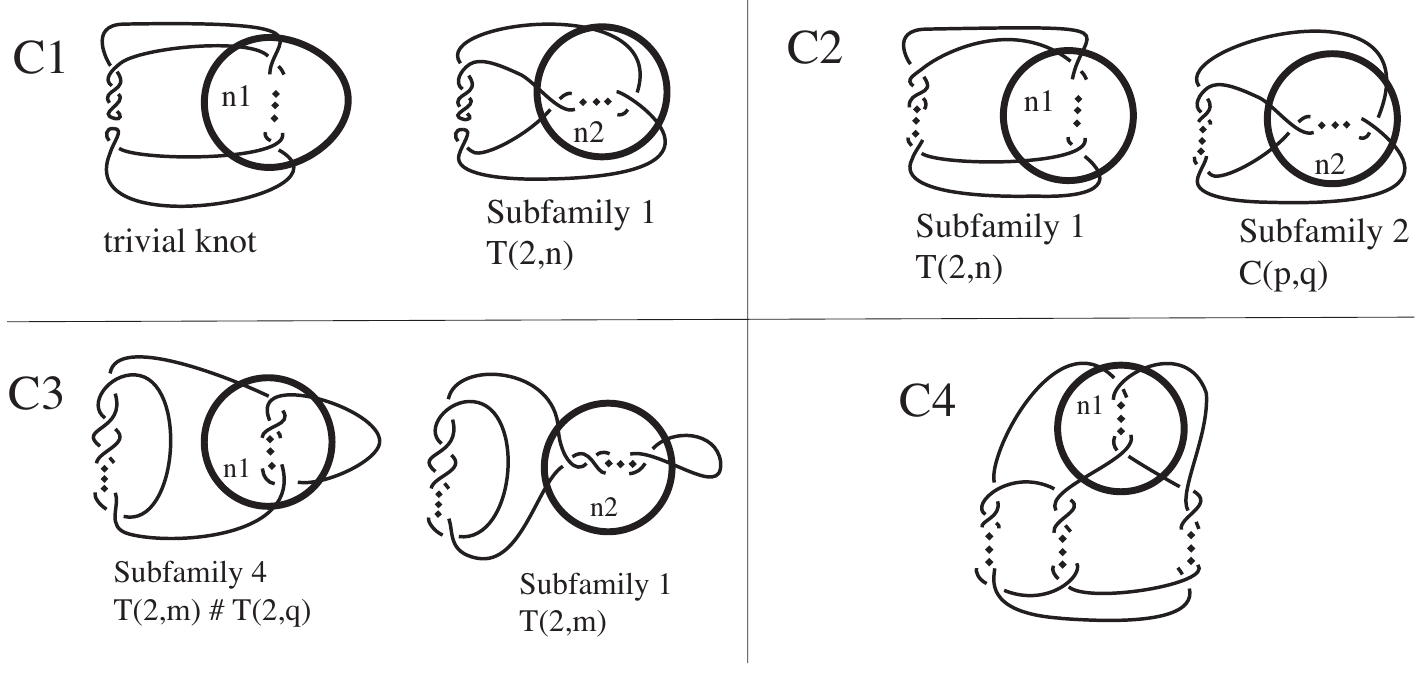}
\caption{Products of recombination with serine recombinases.}
\label{serine}
\end{figure}

\begin{thm}\label{T:serine}
Suppose that Assumptions 1, 2, and 3 hold for a particular serine
recombinase-DNA complex. If the substrate is an
unknot then the only nontrivial products are $T(2,n)$ or
$C(p,q)$. If the substrate is an unlink, then the
only nontrivial product is $T(2,n)$. If the substrate is $T(2,m)$ then all products are in
the family illustrated in Figure \ref{productfamily}.
\end{thm}

\begin{proof}

We saw that as a result of Assumption~3, after $n$ recombination events with serine recombinases, the fixed projection of $B\cap J$ is ambient isotopic fixing $\partial B$ to Form n1 or n2, illustrated in Figure \ref{nBJPaper}.  Also, by Lemma \ref{L:CcapJ}, $C\cap J$ is ambient isotopic fixing $\partial B$ to one of the four forms illustrated in Figure
\ref{newforms}. We obtain the products of serine recombinase from each of the forms of $C\cap J$ illustrated in Figure \ref{newforms} by replacing $B$ with each of Form n1 and Form n2.  The resulting products are illustrated in Figure \ref{serine}.  Note that when $C\cap J$ has Form C4, then $B\cap J$ must have Form B1.  Hence the post-recombinant form of $B\cap J$ must be of Form n1.  Recall again that if $J$ is an unlink, then
$C\cap J$ must have Form C1, and if $J$ is an unknot then $C\cap J$ must have Form C1 or C2.  The theorem follows from Figure \ref{serine}.
\end{proof}

\medskip

Table 1 summarizes the results of Theorems 1 and 2.

\medskip

{\begin{table}
    \begin{tabular}{||l|c|c||} \hline \hline

Recombinase Type & Substrate &   Nontrivial Products \\
         \hline \hline
         Tyrosine & unknot & $T(2,n)$, $C(2,n)$ \\
        \hline
    & unlink & Hopf link$=T(2,2)$\\
        \hline
    & $T(2,m)$ &  Any from Figure \ref{productfamily} \\
        \hline
Serine & unknot & $T(2,n)$, $C(p,q)$\\
        \hline
    & unlink  & $T(2,n)$\\
        \hline
    &  $T(2,m)$ & Any from Figure \ref{productfamily} \\
  \hline

\end{tabular}

\medskip

\caption{Non-trivial products
predicted by our model.}

\end{table}

\section{The minimal crossing number and our model}

\subsection{Our family grows with $n^3$}
 The minimal crossing number of a DNA knot or link can be
determined experimentally using gel electrophoresis. However, there are 1,701,936 knots with minimal crossing number
less than or equal to 16 \cite{HosThis}, and the
number of knots and links with minimal crossing number $n$ grows exponentially as a
function of $n$ \cite{ES3}.  By contrast, we will now prove that the total number of
knots and links in our product family (Figure \ref{productfamily}) grows
linearly with $n^3$.  Note that, while the knots and links in our family have at most four rows containing $p$, $q$, $r$, and $s$ signed crossings respectively, it does not follow that the minimal crossing number of such a knot or link is $|p|+|q|+|r|+|s|$. If the knot or link is not alternating, it is quite possible that the number of crossings can be significantly reduced.  Thus there is no reason to a priori believe that the number of
knots and links in our product family should grow
linearly with $n^3$.

We begin with some results about minimal crossing number which will be used in our proof. We shall denote the {\it minimal crossing
number} of a knot or link $K$ by MCN$(K)$.

\medskip

\begin{lemma} \label{pretzel}   Let $|r|>1$ and $|s|>1$.  Then
$C(r,s)$ is equivalent to both $K(r+1, -1, s+1)$ and $K(r-1, 1,
s-1)$. Furthermore, if $r$ and $s$ have the same sign
then $\mathrm{MCN}(C(r,s))=|r|+|s|-1$, and if $r$ and $s$ have
opposite sign then $\mathrm{MCN}(C(r,s))=|r|+|s|$.
\end{lemma}
\medskip

\begin{proof}
 Figure \ref{C(p,q)}, we show that $C(r,s)$ is ambient isotopic to
$K(r+1, -1, s+1)$ by moving the highlighted strand in front of the
diagram and then turning the horizontal row of $s$ half-twists so that
they become vertical. Analogously, by moving the highlighted
strand behind rather than in front of the rest of the diagram, we
see that $C(r,s)$ is also ambient isotopic to $K(r-1, 1, s-1)$.

\begin{figure}[htpb]
\includegraphics{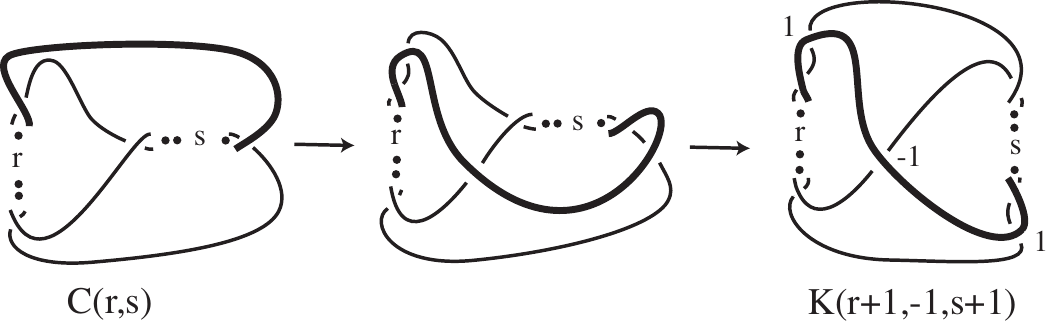}
\caption{An isotopy from $C(r,s)$ to $K(r+1, -1, s+1)$.} \label{C(p,q)}
\end{figure}

We evaluate $\mathrm{MCN}(C(r,s))$ as follows. Murasugi \cite{Mu}
and Thistlethwaite \cite{Th} proved that any reduced alternating
diagram has a minimal number of crossings. Observe that if $r$ and $s$ have opposite
signs, then the standard diagram of $C(r,s)$ is reduced and
alternating. In this case, $\mathrm{MCN}(C(r,s))=|r|+|s|$. If $r$ and $s$ have the same sign and $|r|$, $|s|>1$, then either $r$, $s>1$ or $r$, $s<-1$.  If $r$,
$s>1$, then the diagram of $K(r-1, 1, s-1)$ is reduced and
alternating, since all three rows of crossings are positive. In
this case, $\mathrm{MCN}(C(r,s))=\mathrm{MCN}(K(r-1, 1,
s-1))=r-1+1+s-1=|r|+|s|-1$. If $r$, $s<-1$, then the diagram of
$K(r+1, -1, s+1)$ is reduced and alternating, since all three rows
of crossings are negative. In this case,
$\mathrm{MCN}(C(r,s))=\mathrm{MCN}(K(r+1, -1,
s+1))=-(r+1)+1-(s+1)=|r|+|s|-1$.
\end{proof}
\medskip

To prove our theorem, we will also make use of the following theorem of Lickorish and Thistlethwaite.
\medskip

{\bf  Theorem} \cite{LT}  {\it Suppose that a knot or link $L$ has a
projection as in Figure \ref{Lickorish} with $k\geq 3$, and for each $i$, $R_i\cap L$ is a
reduced alternating projection which contains a crossing between the two arcs at the bottom of $R_i$ (as in Figure \ref{Lickorish}) and at least one other crossing. Then the projection of $L$
has a minimal number of crossings.}

\begin{figure}[htpb]
\includegraphics{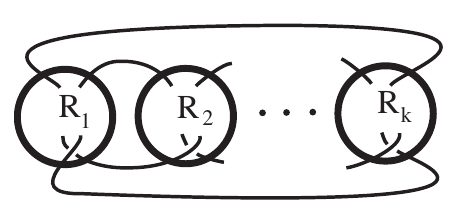}
\caption{Each $R_i$ is reduced, alternating, and has at least two crossings.}
\label{Lickorish}
\end{figure}
\medskip

We shall adopt the language of Lickorish and Thistlethwaite and refer to a projection of the form described by their theorem as a {\it reduced Montesinos diagram}.  Thus by the theorem, any projection of a knot or link which is a reduced Montesinos diagram has a minimal number of crossings.

\begin{thm}\label{n3}
The number of distinct knots and links in the product family illustrated in Figure \ref{productfamily} which have MCN$=n$ grows linearly with $n^3$.
\end{thm}

\begin{proof} We begin by fixing $n$, and suppose that $K$ is a knot or link projection in the family of Figure \ref{productfamily} which has MCN$=n$.  This projection has $|p|+|q|+|r|+|s|$ crossings.  If the given projection of $K$ is reduced alternating or reduced Montesinos, then $|p|+|q|+|r|+|s|=n$.  Otherwise, we show that $K$ is ambient isotopic to one of 24 possible projections  which have a minimal number of crossings.  We will then show that there are at most $96n^3$ possible knots and links in our family with MCN$=n$.

The following example illustrates the type of strand move we shall use to reduce the number of crossings whenever the diagram is neither reduced alternating nor reduced Montesinos.  Observe that the part of our knot or link consisting of the rows containing $r$ and $s$ crossings is alternating if and only if $r$ and $s$ have opposite signs.  If $r$ and $s$ have the same sign, then by moving a single strand  (as in Figure \ref{rs}), this part of the knot or link becomes alternating.  This isotopy removes a crossing from both the $r$ row and the $s$ row and adds a single new crossing.  Thus we reduce this part of the diagram from having $|r|+|s|$ crossings in a non-alternating form to having $(|r|-1)+(|s|-1) +1$ crossings in an alternating form.   All of the isotopies we use to to get rid of unnecessary crossings involve moving at most three such strands.

\begin{figure}[htpb]
\includegraphics{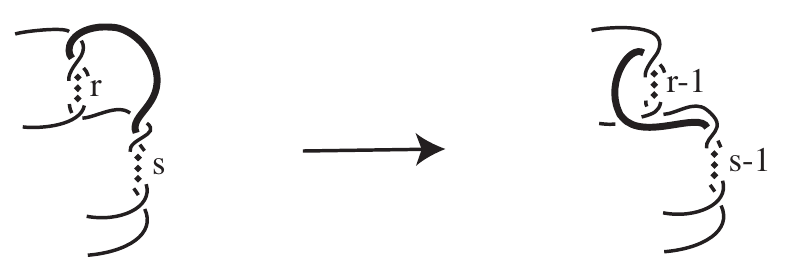}
\caption{By moving a single strand we  reduce from $|r|+|s|$ crossings originally to $(|r|-1)+(|s|-1) +1$ crossings in the alternating diagram.}
\label{rs}
\end{figure}

Next we will discuss the one exceptional case where we cannot obtain a reduced alternating or reduced Montesinos diagram by moving some strands of $K$.  This is the case when $K$ is a knot or link in our family with $r>1$, $p$, $q<-2$, and $s=1$. In its original form, the projection has $-p-q+r+1$ crossings.  We can move a single strand of the diagram to obtain a projection with only $-p+(-q-1)+(r-1)+1$ crossings (illustrated on the left in Figure \ref{exception}).  We define a {\it Hara-Yamamoto} projection as one in which there is a row of at least two crossings which has the property that if this row is cut off from the rest of the projection and the endpoints are resealed in the two natural ways, then both resulting projections are reduced alternating.  The projection on the left of Figure \ref{exception} is Hara-Yamamoto because the projections (on the right) obtained by resealing the endpoints are both reduced alternating.  Hara and Yamamoto \cite{HY} have shown that any Hara-Yamamoto projection has a minimum number of crossings.  Thus the projection on the left of Figure \ref{exception} has a minimal number of crossings.

\begin{figure}[htpb]
\includegraphics{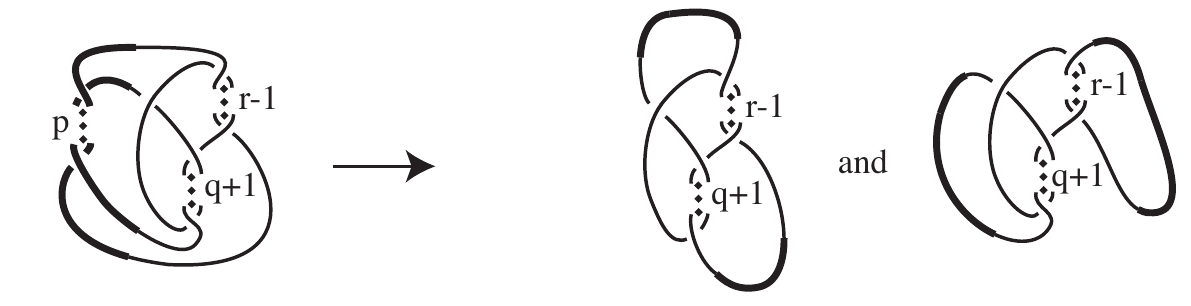}
\caption{If we cut off the row of $p$ crossings on the left and reseal the strands in the two natural ways, then both resulting projections are reduced alternating.}
\label{exception}
\end{figure}

We will consider 27 cases according to the values of $p$, $q$, $r$, and $s$, and show that in all but the above exceptional case $K$ is isotopic to a diagram that is either reduced alternating or reduced Montesinos and hence has minimal crossing number. Since there are so many cases, we display the results in a table rather than discussing each case individually.  We make the following notes about the table. In the second column we list the form of the knot or link which has a minimal number of crossings (e.g. reduced alternating). If the knot or link is isotopic to a clasp, pretzel, or torus knot or link we will list the specific form (e.g. $T(2,n)$).  If the minimal crossing form is either a clasp $C(r,s)$ or a pretzel of the form $K(r,\pm1,s)$ then (according to  Lemma \ref{pretzel}) which one of these is the minimal crossing form depends on the signs of $r$ and $s$.  In this case, we just list one of these two forms though the one we list is not necessarily the form with the fewest  number of crossings, as we do not know the signs or specific values of the variables.  In this case, for the MCN we write an expression with {\bf (-1?)} at the end to mean that depending on the relevant variables the MCN may be one smaller.  If one of these knots or links contains a trivial component, we use the shorthand {\it $+$ O} to indicate this in the table.   We shall consider a knot or link and its mirror image to be of the same link type, and hence we will not count both.  Thus without loss of generality, we shall assume that $r\geq 0$.  Also, observe that the rows of crossings containing $p$ and $q$ crossings are intechangeable in Figure \ref{productfamily}, so we treat the variables $p$ and $q$ as interchangeable.  We list the MCN as an unsimplified function of $p$, $q$, $r$, and $s$ to help the reader recreate the isotopy taking the original form to the minimal crossing form. Finally, apart from the cases where $K$ reduces to $T(2,m)$ or $C(2,m)$, we obtain the upper bounds for the number of links in each case by expressing MCN$=n$ as a sum of nonnegative integers.  This enables us to find an upper bound for the number of knots and links with MCN$=n$ in each case.  Note that the upper bounds given are intended to be simple rather than as small as possible. In particular, a number of our cases overlap, and thus some knots and links are counted more than once.
Also, for certain specific values of $p$, $q$, $r$, and $s$, we may obtain a trivial knot or link.  However, we do not specifically exclude these cases from our table.

\begin{footnotesize}
\begin{table}[h]
    \begin{tabular}{{||l||l|c|l|c||}} \hline \hline
Values of $p$, $q$, $r$, $s$ &  Minimal crossing form & Strands &MCN written as a sum  & Upper bound \\ for $r\geq 0$ && moved &  of nonnegative integers& on \# of links\\
         \hline \hline
         $p=q=0$ &  $C(r,s)+$O  & 0 & $r+|s|$ \hfill{\bf (-1?)}  & $4n$\\
          \hline
         $r=0$ &  $T(2,p+q)$ & 0& $|p+q|$ & 1 \\
      \hline

           $r=1$, $p\not=0$, $q=0$ &  $T(2,p)\#T(2,s+1)$ & 0 & $|p|+|s+1|$ & $2n$ \\
           \hline
            $r=1$, $p\not =0$, $q\not =0$ & $K(p,q,s+1)$  & 0 & $|p|+|q|+|s+1|$ \hfill{\bf (-1?)}& $8n^2$  \\
           \hline
$r>1$, $p\not=0$, $q=0$ &  $T(2,p)\#C(r,s)$  & 0 & $|p|+r-s$ \hfill{\bf (-1?)}& $8n^2$ \\
           \hline

           $r>1$, $pq=-1$ & $T(2,r)$ & 0& $r$ & 1\\
          \hline
          $r>1$, $pq=1$, $s=0$ & $C(\pm 2, r)$  & 0 & $2+r$ \hfill{\bf (-1?)} & 2\\
          \hline

          $r>1$, $p \geq 1$, $q=1$, $s>0$ & reduced alternating & 1 & $p+ (r-1)+(s-1)+2$ & $n^2$ \\
          \hline
          $r>1$, $p=q=1$, $s<0$ & reduced alternating &1& $r+(-s-1)+2$&$n$ \\
          \hline

          $r>1$, $p\leq -1$, $q=-1$, $s>0$ & reduced alternating & 2 & $-p+(r-1)+(s-2)+2$& $n^2$\\
          \hline
            $r>1$, $p$, $q<0$, $s\leq 0$ & reduced alternating & 0&$-p-q+r-s$ &$n^3$\\
          \hline
          $r>1$,  $p$, $q>1$, $s=0$ & reduced alternating & 2& $(p-1)+(q-1)+(r-2)+2$&$n^2$ \\
          \hline

          $r>1$, $p<-1$, $q>1$, $s=0$ & reduced alternating & 1& $-p+(q-1)+(r-1)+1$ &$n^2$\\
          \hline

           $r>1$, $|p|>1$, $|q|=1$, $s=0$ & $C(r\pm 1,p)$  & 0& $-p+(r\pm1)$ \hfill{\bf (-1?)} &$4n$ \\
          \hline
          $r>1$, $qs=-1$ & $T(2, r+p\pm1)$ & 1& $|r+p\pm1|$ & 1 \\
          \hline
          $r>1$, $p>0$, $q=1$, $s<0$ & reduced alternating  & 1& $p+r+(-s-1)+1$& $n^2$\\
          \hline

          $r>1$, $p\leq-2$, $q=1$, $s\leq-2$ & reduced alternating & 1& $(-p-2)+r+(-s-2)+1$ & $n^2$\\
          \hline

             $r>1$, $p$, $q>0$, $s=1$ & reduced alternating & 1& $p+q+(r-1)+1$ &$n^2$\\
          \hline

                      $r>1$, $p<-1$, $q>0$, $s=1$ & reduced alternating & 1& $(-p-1)+q+(r-1)$ & $n^2$ \\
          \hline

          $r>1$, $p<-1$, $q=1$, $s>1$ & reduced alternating & 2& $(-p-1)+(r-1)+(s-1)+2$ &$n^2$\\
          \hline

           $r>1$, $p>1$, $q=-1$, $s<0$ & reduced alternating  & 1& $(p-1)+r-s+1$& $n^2$\\
          \hline
           $r>1$, $p>1$, $q=-1$, $s=2$ & trivial  & 2 & $0\not =n$& 0\\
          \hline

$r>1$, $p>1$, $q=-1$, $s>2$ & reduced alternating  & 3 & $(p-2)+(r-1)+(s-3)+2$& $n^2$\\
          \hline

                $r>1$, $|p|$, $|q|>1$, $s<0$ &  reduced Montesinos & 0& $|p|+|q|+r-s$ &$4n^3$\\
          \hline
          $r>1$, $|p|$, $|q|>1$, $s>1$ &  reduced Montesinos & 1& $|p|+|q|+(r-1)+(s-1)+1$ &$4n^3$\\
          \hline

       $r>1$, $p<-1$, $q=-2$, $s=1$ & $K(p,2,r-1)$  & 1& $-p+2+(r-1)$ & $n$   \\
       \hline

       $r>1$, $p$, $q<-2$, $s=1$ & Hara-Yamamoto & 1& $-p+(-q-1)+(r-1)$ & $n^2$\\
       \hline

        \hline
        \hline

\end{tabular}

\medskip

\caption{The minimal crossing forms of links in the family in Figure \ref{productfamily}}
\end{table}
\end{footnotesize}

There are 26 non-trivial cases in the table.  However, all three instances of a $T(2,m)$ yield the same knot or link.  Thus there are at most 24 distinct families of knots and links listed in the table. The number of knots and link in each of these families is bounded above by $4n^3$ (in fact, for most of the cases there are significantly fewer than $4n^3$ knot and link types).  It follows that for a given $n$, the number of distinct knots and links in our product family (Figure \ref{productfamily} which have MCN$=n$ is bounded above by $24\times 4n^3=96n^3$.  In particular, the number of distinct knots and links with the form of Figure \ref{productfamily} which have MCN$=n$ grows linearly with $n^3$.

\end{proof}

\medskip

It follows from Theorem \ref{n3} that the proportion of all knots and links
which are contained in our family decreases exponentially as $n$
increases. Thus, for a knotted or linked product, knowing the MCN
and that it is constrained to this family allows us to
significantly narrow the possibilities for its precise knot or
link type. The model described herein thus provides an important
step in characterizing DNA knots and links which arise as products
of site-specific recombination.
\medskip

\subsection{Products whose MCN is one more than the substrate}
Finally, we prove a more directly applicable theorem as follows.  Site-specific recombination often adds a single crossing to the MCN of a knotted or linked substrate.  If the substrate is $T(2,m)$ and the product of a single recombination event has $\mathrm{MCN}=m+1$, then we can further restrict the resulting knot or link type.
 \medskip

\begin{thm}\label{T:MCN}
Suppose that Assumptions 1, 2, and 3 hold for a particular
recombinase-DNA complex with substrate $T(2,m)$,
with $m>0$. Let $L$ be the product of a single recombination
event, and suppose that $\mathrm{MCN}(L)=m+1$. Then $L$ is either
$T(2,m+1)$, $C(-2,m-1)$, or $K(s,t,1)$ with $s,t>0$ and $s+t=m$ (see Figure  \ref{mcntheorem}).
\end{thm}

  \begin{figure}[htpb]
\includegraphics{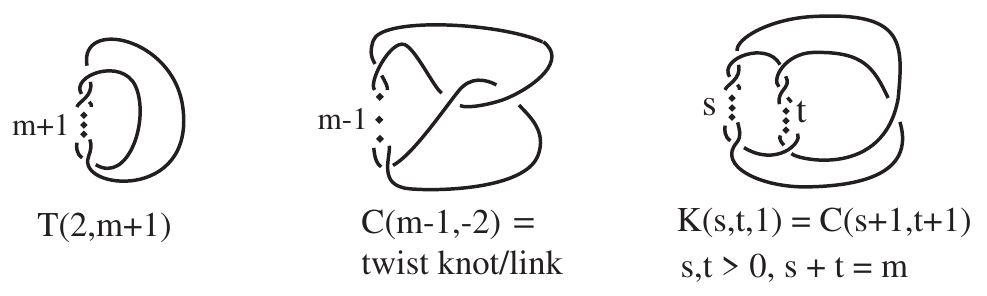}
\caption{These are the only possible products, if the substrate is $T(2,m)$ and the product has MCN$=m+1$.}
\label{mcntheorem}
\end{figure}

\begin{proof}  For $m=1$, $T(2,1)$ is the trefoil knot $3_1$, and hence  $L$ must be the figure eight knot $4_1$ which can also be written as $K(2,1,1)$.   Thus from now on we assume that $m\geq 2$.  By Assumption~1, there a projection of $J$ such that $B\cap J$ has at most one crossing.  Since $J=T(2,m)$, the proof of Lemma \ref{L:CcapJ} shows that $C\cap J$ is ambient isotopic, fixing $\partial B$ to a projection with Form C2, C3, or C4 (see
Figure \ref{newforms}). Furthermore, when $C\cap J$ has Form C4, then $p+q=m$.  By Assumption~3, after a single recombination event with either serine or tyrosine recombinases the post-recombinant form of $B\cap J$ is ambient isotopic fixing $\partial B$ to one of  those illustrated in Figure \ref{tyroBJ1}. Thus any knotted or linked product $L$ has one of the forms illustrated in
Figure \ref{tyrosine}.

First suppose that $L$ has one of the forms illustrated when $C\cap J$ has Form C2 or C3.  We see that $L$ cannot be $T(2,m)\#T(2,2)$, since MCN$(T(2,m)\#T(2,2))=m+2$.  Certainly, $L$ cannot be $T(2,m)$ with a trivial component.  If $L=T(2,n)$ then $n=m+1$, so we are done.  If $L=C(-2,n)$ and $n>1$, then $n=m-1$, so again we are done.  If $L=C(2,n)$ and $n>1$, then $L=K(1,1,n-1)$. In this case $n=m$, and again we are done.

Now suppose that $L$ has one of the forms illustrated when $C\cap J$ has Form C4. If $L=K(p,q,a)$ for some value of $a$, then $L$ has a projection in the product family illustrated in Figure \ref{productfamily} with $r=1$.  Otherwise, $L$ is a member of the product family with $r=\pm2$.  However, if $r=-2$, we can turn over the top loop to get $r=2$ (this will also add on positive crossing to the $s$ row).  Thus we shall now assume that $L$ has a projection in the product family (i.e., Figure \ref{productfamily}) with either $r=1$ or $r=2$. Table 2 lists all of the nontrivial knots and links in this family when $r\geq 0$.  Thus all of the products that we are considering occur in Table 2.  We would like to know which of the cases in Table 2 have $r=1$ or $2$, $p+q\geq 2$, and MCN$=p+q+1$.  The following table answers this question.

\begin{footnotesize}
\begin{table}[h]
    \begin{tabular}{{||l|l|l|c|c||}} \hline \hline
Values of $p$, $q$, $r$, $s$ &  Minimal crossing form & MCN written as a sum&Is $r=1$ or $2$  & Can MCN$=$ \\ for $r\geq 0$ & &  of nonnegative integers& and $p+q\geq 2$?& $p+q+1$?\\
         \hline \hline
         $p=q=0$ &  $C(r,s)+$O   & $r+|s|$  \hfill{\bf (-1?)} & no & -\\
          \hline
         $r=0$ &  $T(2,p+q)$ & $|p+q|$ & no & - \\
      \hline

           $r=1$, $p\not=0$, $q=0$ &  $T(2,p)\#T(2,s+1)$  & $|p|+|s+1|$ & yes & no\\
           \hline
            $r=1$, $p\not =0$, $q\not =0$ & $K(p,q,s+1)$   & $|p|+|q|+|s+1|$ \hfill{\bf (-1?)} & yes & if $s=0$ \\
           \hline
$r>1$, $p\not=0$, $q=0$ &  $T(2,p)\#C(r,s)$   & $|p|+r+|s|$ \hfill{\bf (-1?)}& yes &no\\
           \hline

           $r>1$, $pq=-1$ & $T(2,r)$ & $r$ & no& -\\
          \hline
          $r>1$, $pq=1$, $s=0$ & $C(\pm 2, r)$   & $2+r$ \hfill{\bf (-1?)}& yes & if $r=p+q-1$\\
          \hline

          $r>1$, $p \geq 1$, $q=1$, $s>0$ & reduced alternating  & $p+ (r-1)+(s-1)+2$ & yes &no\\
          \hline
          $r>1$, $p=q=1$, $s<0$ & reduced alternating & $r+(-s-1)+2$& yes & no\\
          \hline

          $r>1$, $p\leq -1$, $q=-1$, $s>0$ & reduced alternating  & $-p+(r-1)+(s-2)+2$& no& -\\
          \hline
            $r>1$, $p$, $q<0$, $s\leq 0$ & reduced alternating &$-p-q+r-s$ & no& -\\
          \hline
          $r>1$,  $p$, $q>1$, $s=0$ & reduced alternating & $(p-1)+(q-1)+(r-2)+2$& yes & only if $r=3$\\
          \hline

          $r>1$, $p<-1$, $q>1$, $s=0$ & reduced alternating & $-p+(q-1)+(r-1)+1$ & yes& no\\
          \hline

           $r>1$, $|p|>1$, $|q|=1$, $s=0$ & $C(r\pm 1,p)$  & $-p+(r\pm1)$ \hfill{\bf (-1?)}& yes & yes\\
          \hline
          $r>1$, $qs=-1$ & $T(2, r+p\pm1)$ & $|r+p\pm1|$  & yes& yes\\
          \hline
          $r>1$, $p>0$, $q=1$, $s<0$ & reduced alternating  & $p+r+(-s-1)+1$& yes& no\\
          \hline

          $r>1$, $p\leq-2$, $q=1$, $s\leq-2$ & reduced alternating & $(-p-2)+r+(-s-2)+1$ &  no& -\\
          \hline

             $r>1$, $p$, $q>0$, $s=1$ & reduced alternating & $p+q+(r-1)+1$ & yes& no\\
          \hline

                      $r>1$, $p<-1$, $q>0$, $s=1$ & reduced alternating & $(-p-1)+q+(r-1)$ & yes & no\\
          \hline

          $r>1$, $p<-1$, $q=1$, $s>1$ & reduced alternating & $(-p-1)+(r-1)+(s-1)+2$ & no& -\\
          \hline

           $r>1$, $p>1$, $q=-1$, $s<0$ & reduced alternating  & $(p-1)+r-s+1$& yes& no\\
          \hline
           $r>1$, $p>1$, $q=-1$, $s=2$ & trivial   & $0\not =n$& yes& no\\
          \hline

$r>1$, $p>1$, $q=-1$, $s>2$ & reduced alternating   & $(p-2)+(r-1)+(s-3)+2$& yes&no\\
          \hline

                $r>1$, $|p|$, $|q|>1$, $s<0$ &  reduced Montesinos & $|p|+|q|+r-s$ & yes&no\\
          \hline
          $r>1$, $|p|$, $|q|>1$, $s>1$ &  reduced Montesinos & $|p|+|q|+(r-1)+(s-1)+1$ & yes & no\\
          \hline

       $r>1$, $p<-1$, $q=-2$, $s=1$ & $K(p,2,r-1)$  & $-p+2+(r-1)$ & no   & -\\
       \hline

       $r>1$, $p$, $q<-2$, $s=1$ & Hara-Yamamoto & $-p+(-q-1)+(r-1)$ & no& -\\
       \hline

        \hline
        \hline

\end{tabular}

\medskip

\caption{Which products can have MCN$=p+q+1$?}
\end{table}
\end{footnotesize}
\medskip

The only subtle case in the table is where $L=C(r\pm1,p)$. In this case we must have  $r=2$, $|q|=1$, and $p+q\geq 2$.  It follows that $p\geq 1$. Since $L$ is non-trivial, we must have $L=C(3,p)=K(2,1,p-1)$.  Now MCN$(L)=m+1$ implies that $p-1+2=m=p+q$.  Thus $L=K(s,t,1)$ where $s+t=p+q$.  Now from Table 3, we can see that if $r=1$ or $2$, $p+q\geq 2$, and MCN$(L)=p+q+1$, then $L$ is either $T(2,m+1)$, $C(m-1,-2)$, or $K(s,t,1)$ with
$s$, $q>0$, and $s+t=m$.
\end{proof}
\medskip

We now illustrate an application of Theorem~3 (further applications of our model are discussed in \cite{BFbio}). Bath \textit{et
al} used the links $6^2_1$ and $8^2_1$ as the substrates for Xer
recombination, yielding a knot with MCN=7 and a knot with MCN=9,
respectively. These products have not been characterized beyond
their minimal crossing number, and MCN is not sufficient to
determine the knot type. In particular, there are seven knots with
MCN=7 and 49 knots with MCN=9.

Theorem \ref{T:MCN} significantly reduces the number of
possibilities for each of these products. In particular, it
follows from Theorem \ref{T:MCN} that the 7-crossing products of
Xer must be $7_1=T(2,7)$, $7_2=C(5,-2)$ or $7_4=K(3,3,1)$; and the
9-crossing products of Xer must be $9_1=T(2,9)$, $9_{2}=C(7,-2)$,
or $9_5=K(5,3,1)$. All of these possible products are actually
4-plats. This example shows how our model complements the work of
\cite{Dar}, which restricts attention to the tangle model, and
thus assume that all products are 4-plats. In \cite{new}, building
on earlier work of \cite{Vaz2,Dar,us,us2}, we use our model
together with tangle calculus to completely classify all tangle
solutions to these Int-Xer equations.

\medskip

\section{Acknowledgements}

The authors wish to thank Andrzej Stasiak, De
Witt Sumners, Alex Vologodskii, and Stu Whittington for helpful
conversations. Dorothy Buck was partially supported by Grant \# DMS-0102057
from the National Science Foundation's Division of Mathematical
Sciences. Erica Flapan was partially supported by an Association for Women
in Mathematics Michler Collaborative Research Grant.


\bibliographystyle{amsplain}

\end{document}